\documentclass[12pt,a4paper]{article}

\usepackage{amssymb,amsmath,amsthm}
\usepackage[english]{babel}
\usepackage{t1enc}
\usepackage[latin2]{inputenc}
\usepackage{epsfig}
\usepackage{subfigure}
\usepackage{footnpag}

\newtheorem{thm}{Theorem}

\newtheorem{claim}[thm]{Claim}

\newtheorem*{problem}{Problem}

\usepackage{indentfirst}
\def\R{\mathbb{R}}

\begin{document}

\title{Octants are Cover-Decomposable into Many Coverings}
\author{Bal\'azs Keszegh \thanks{Research supported by OTKA under grant NN 102029 (EUROGIGA project GraDR 10-EuroGIGA-OP-003) and under grant NK 78439.}\\
\and D\"om\"ot\"or P\'alv\"olgyi\thanks{Research supported by Hungarian National Science Fund (OTKA), under grant PD 83586 and under grant NN 102029 (EUROGIGA project GraDR 10-EuroGIGA-OP-003) and the J\'anos Bolyai Research Scholarship of the Hungarian Academy of Sciences.}
}


\maketitle

\begin{abstract}
We prove that octants are cover-decomposable into multiple coverings, i.e., for any $k$ there is an $m(k)$ such that any $m(k)$-fold covering of any subset of the space with a finite number of translates of a given octant can be decomposed into $k$ coverings. As a corollary, we obtain that any $m(k)$-fold covering of any subset of the plane with a finite number of homothetic copies of a given triangle can be decomposed into $k$ coverings. 
Previously only some weaker bounds were known for related problems \cite{V10}.
\end{abstract}

\medskip

\section{Introduction}
Let ${\cal P}=\{\ P_i\ |\ i\in I\ \}$ be a collection of geometric sets in $\R^d$.
We say that ${\cal P}$ is an {\em $m$-fold covering} of a set $S$ if every point of $S$ is contained in at least $m$ members of $\cal P$. A $1$-fold covering is simply called a {\em covering}. 

\medskip

\noindent {\bf Definition.} A geometric set $P\subset \R^d$ is said to be {\em cover-decomposable into $k$ coverings} if there exists a (minimal) constant $m(k)=m_P(k)$ such that every $m(k)$-fold covering of any $X\subset\R^d$ with a finite number of translates of $P$ can be decomposed into $k$ coverings of $X$.
If $m=m(2)$ exists, we say that $P$ is  {\em cover-decomposable}.
If $m(k)$ exists for every $k$, we say that $P$ is {\em cover-decomposable into many coverings}.

\medskip

We note that in the literature (\cite{P10}, \cite{PPT11}) the definition is slightly different and the notion defined here is sometimes called {\em finite-cover-decomposable}, however, to avoid unnecessary complications, we simply call it {\em cover-decomposable}.

The main result of this paper is

\begin{thm}\label{newmainthm} Octants are cover-decomposable into many coverings, any $m(k)$-fold covering of any subset of $\R^3$ with a finite number of translates of a given octant can be decomposed into $k$ coverings, where $m(k)\le 12^{2^{k-1}+2^{k-2}-2}<12^{2^k}$.
\end{thm}

We prove this theorem in Section 2. In the remainder of this section we give a brief history of the problem. For a more detailed introduction and other results on cover-de\-com\-po\-sa\-bi\-li\-ty, see the recent surveys \cite{P10} and \cite{PPT11} and the papers \cite{A08,B07,GV10,P80,P86,PT07,P11,PT10,TT07}.

It is easy to see that a quadrant (i.e., a $2$-dimensional orthant) is cover-decomposable. Cardinal \cite{C11} noticed that orthants in $4$ and higher dimensions are not cover-decomposable as there is a plane on which their trace can be any family of axis-parallel rectangles and it was shown by Pach, Tardos and T\'oth \cite{PTT09} that such families might not be decomposable into two coverings. Cardinal asked whether octants ($3$-dimensional orthants) are cover-decomposable and this was settled in \cite{KP} by the following theorem.
 
\begin{thm}\label{mainthm} Octants are cover-decomposable, any $12$-fold covering of any subset of $\R^3$ with a finite number of translates of a given octant can be decomposed into two coverings.
\end{thm}

In fact, the following equivalent, dual form of Theorem \ref{mainthm} was proved.

\begin{thm}\label{dualthm} Any finite set of points in $\R^3$ can be colored with two colors such that any translate of a given octant with at least $12$ points contains both colors.
\end{thm}

Similarly, we will prove the following, equivalent form of Theorem \ref{newmainthm}, using Theorem \ref{dualthm}.

\begin{thm}\label{newdualthm} Any finite set of points in $\R^3$ can be colored with $k$ colors such that any translate of a given octant with at least $12^{2^{k-1}+2^{k-2}-2}$ points contains all $k$ colors.
\end{thm}

\section{Proof of Theorem 1} \label{secmain}
In the next section we give the basic definitions and notations already established in our previous paper \cite{KP}, so the reader familiar with these results can skip straight to Section \ref{theproof}.

\subsection{Definitions and Notations}
Denote by $W$ the octant with apex at the origin containing $(-\infty,-\infty,-\infty)$.
We will work in the dual setting, that is we have a finite set of points, $P$,
in the space, that we want to color with $k$ colors such that any translate of $W$ with at least $m(k)$ points contains all $k$ colors. We call such a $k$-coloring of a point set in the space a {\em $k$-good coloring}. If such a coloring exists for any $P$, then it follows using a standard dualization argument (see \cite{P10} or \cite{PPT11}) that $W$ (and thus any octant) is cover-decomposable into $k$ coverings. So from now on our goal will be to show the existence of such a coloring.

For simplicity, suppose that no number occurs multiple times among the coordinates of the points of $P$ (otherwise, by a small perturbation of $P$ we can get such a point set, and its coloring will be $k$-good for $P$). Denote the point of $P$ with the $t^{th}$ smallest $z$ coordinate by $p_t$ and the union of $p_1,\ldots,p_t$ by $P_t$. First we will show how to reduce the coloring of $P$ to a planar and thus more tractable problem.

Denote the projection of $P$ on the $z=0$ plane 
by $P'$. Similarly denote the projection of $p_t$ by $p_t'$, the projection of $P_t$ by $P_t'$ and the projection of $W$ by $W'$. Therefore $W'$ is the quadrant with apex at the origin containing $(-\infty,-\infty)$.

For such an ordered planar point set $P'$ we say that a coloring with $k$ colors of it is a {\em $k$-good coloring}, if for any $t$ and any translate of $W'$ containing at least $m(k)$ points of $P_t'$, it is true that the intersection of this translate and $P_t'$ contains all $k$ colors. We use the same notation ($k$-good) for two differently defined colorings, because a $k$-good coloring of a spatial point set and a $k$-good coloring of the corresponding planar point set are equivalent problems.

\begin{claim}\label{project} The ordered planar point set $P'$ has a $k$-good coloring if and only if the spatial point set $P$ has a $k$-good coloring.
\end{claim}
 
We omit the proof, as it is not too hard and it is a straightforward generalization of the respective Claim from \cite{KP}. 
 
Now we will prove that any $P'$ has a $k$-good coloring, thus establishing Theorem \ref{newdualthm} and since they are equivalent, also Theorem \ref{newmainthm}.
To avoid going mad, we will omit the apostrophe in the following, so we will simply write $W$ instead of $W'$ and so on.
Also, we will use the term {\em wedge} to denote a translate of $W$. 

A possible way to imagine this planar problem is that in every step $t$ we have a set of points, $P_t$, and our goal is to color the coming new point, $p_{t+1}$, such that we always have a $k$-good coloring.
We note that this would be impossible in an online setting, i.e. without knowing in advance which points will come in which order.
(For related problems, see \cite{KNP}.)
But using that we know in advance every $p_i$ makes the problem solvable.

We introduce some notation.
If $p_x<q_x$ but $p_y>q_y$ then we
say that $p$ is NW from $q$ and $q$ is SE from $p$. In this case we call $p$ and $q$
incomparable. Similarly, $p$ is SW from $q$ (and $q$ is NE from $p$) if and only if 
both coordinates of $p$ are smaller than the respective coordinates of $q$.



\subsection{The Proof}\label{theproof}
Here we prove by induction on $k$ that any ordered planar point set $P$ has a $k$-good coloring.
For $k=2$ it follows from Theorem \ref{dualthm} that $m(2)\le 12$.
Now we will prove that $m(k)\le 144((m(k-1))^2-m(k-1))+1$, establishing Theorem \ref{newdualthm}.
We will start with a preprocessing part where we give an algorithm that partitions the point set, then using the partition we can define the coloring algorithm. Finally we will prove that it indeed gives a $k$-good coloring.

\begin{figure}[t]
    \centering
    \subfigure[The partition to regions $R_{p_i}$.]{\label{fig_regions}
        \includegraphics[scale=0.55]{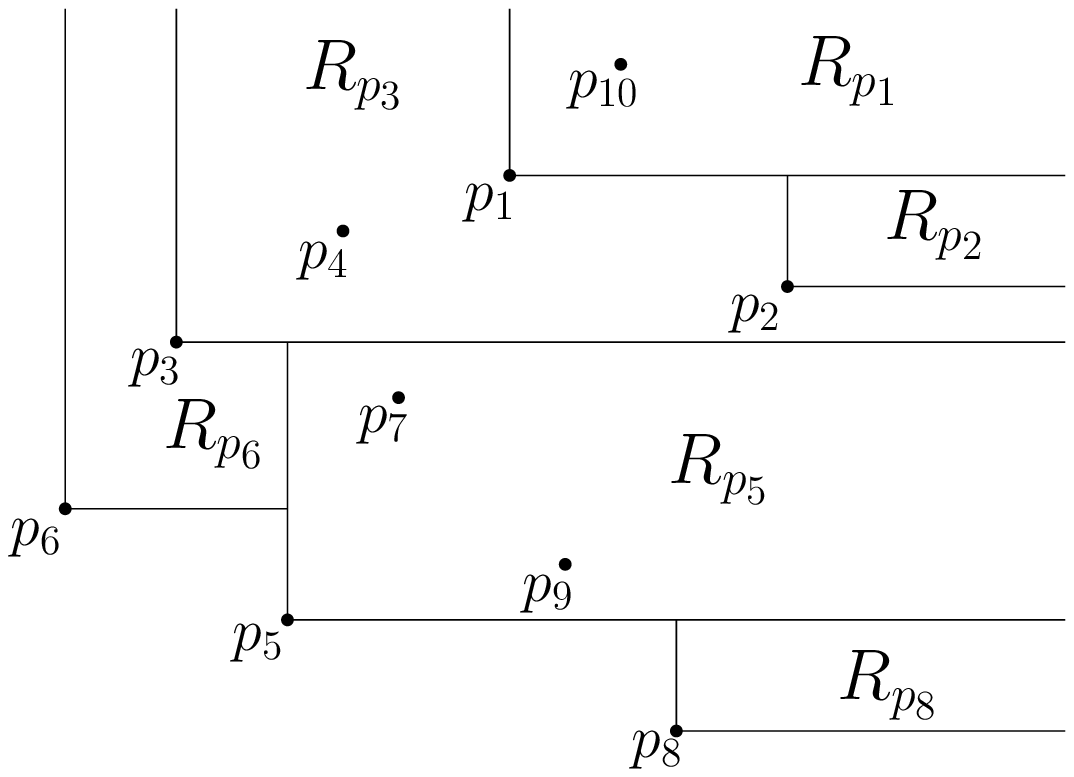}        
        }        
    \hskip 10mm
    \subfigure[Final step of the proof, with $V$ containing many $J$-points.]{\label{fig_important}
        \includegraphics[scale=0.6]{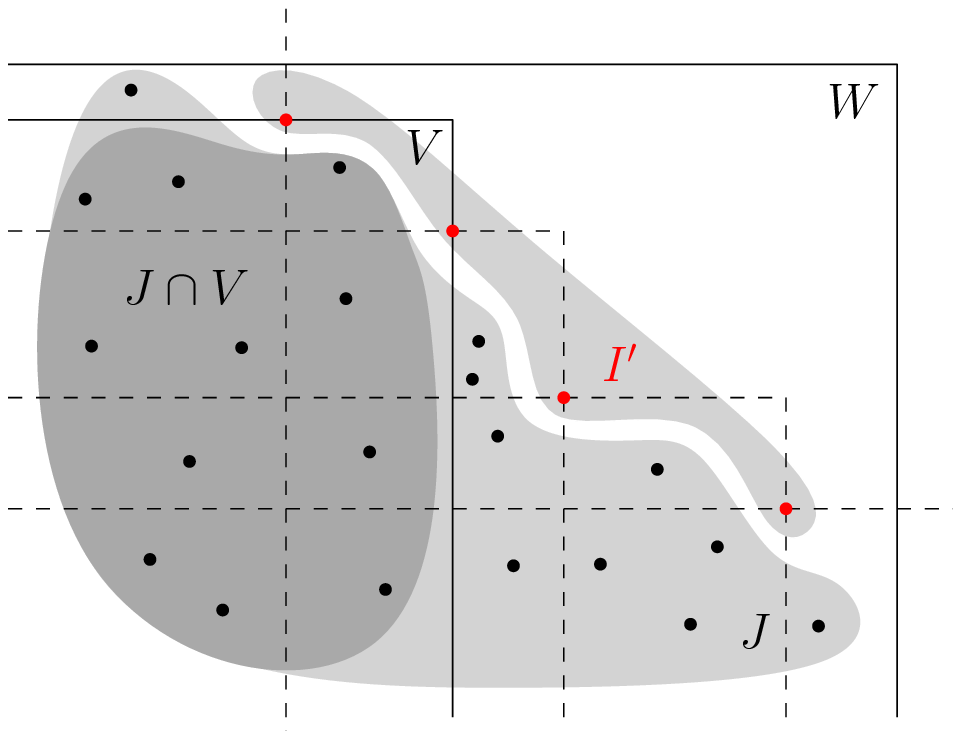}           
        } 
\end{figure}

To distinguish the points of $P$ from other points of the plane, from now on we call them {\em $P$-points}. We similarly use this notation for other sets, e.g. the points of $X$ are called {\em $X$-points}.
First we define an algorithm which partitions the $P$-points into subsets. One set is the set of important points, $S$, and for any $S$-point $p$ we associate a region, $R_p$, that contains only points that are NE from $p$ (but not necessarily all such points) and we denote the $S$-points from $R_p$ by $S_p=P\cap R_p$. Note that all $S_p$-points lie NE from $p$. The regions $R_p$ for $p\in S$ partition those points of the plane that lie NE from some $P$-point. The set $S$ and the sets $S_p$ together partition $P$. Note that as we supposed that no number occurs multiple times among the coordinates of the $P$-points, when we say that the regions partition some other region, we do not need to care about their boundaries. For an illustration see Figure \ref{fig_regions}.

The partitioning algorithm is the following. We process the points one-by-one according to their order in $P$, starting with $p_1$ and recursively define the partition and the regions while we maintain the above properties at all times. We begin by putting $p_1$ into $S$ and set $R_{p_1}$ to be all the points of the plane that are NE from $p_1$, $S_{p_1}$ is empty. Now suppose we are at time $t$ and we process the next point $q=p_t$. If $q$ is in one of the regions $R_p$ (where $p\in P_{t-1}$) then we put $q$ into $S_p$. Otherwise, $q$ is not NE from any of the $P$-points, in which case we put $q$ in $S$. We associate the region $R_q$ with $q$ which contains those points of the plane which lie NE from $q$ but are not contained in any of the regions $R_p$ for $p\in P_{t-1}$ and thus $S_q$ is empty. It is easy to see that all properties are maintained.

Now we can define the coloring algorithm. The set $S$ is first colored by induction using the $k-1$ colors $\{1,2,\dots,k-2,red\}$. Then we recolor the points that are red by induction using two colors, $k-1$ and $k$ (during this recoloring we completely ignore all non-red points). Finally, for each $S_p$ we color $S_p$ by induction using $k-1$ colors, $\{1,2,\dots,k\}\setminus \{color(p)\}$, i.e. with colors that differ from the color of $p$ (again, when coloring some $S_p$, we completely ignore points not in $S_p$).

Now we prove that the algorithm is correct.
Suppose at any time we have a wedge $W$ with $144((m(k-1))^2-m(k-1))+1$ points in it.
If $W$ contains at least $m(k-1)$ $S_p$-points for some $p$, then we are done, because we have $k-1$ different colors in the region, plus the point associated with this region, that has the missing $k^{th}$ color.

Otherwise, since if $W$ contains a point from a region $R_p$, it also contains the $S$-point $p$, we can conclude that $W$ contains at least $\lceil \frac{144((m(k-1))^2-m(k-1))+1}{m(k-1)} \rceil$ $=$ $144m(k-1)-143$ $S$-points.
Since $144m(k-1)-143$ $\ge$ $m(k-1)$, we know from the correctness of the algorithm for $k-1$ colors that $W$ contains all colors that are at most $k-2$.
If $W$ also contains $m(2)=12$ points that were colored red, then after the recoloring it contains points with color $k-1$ and $k$ as well (using the correctness of the algorithm for two colors).
Suppose that this is not the case and denote  by $I$ the at most $11$ $S$-points from $W$ that were first colored red (and then $k-1$ or $k$).
Using the pigeonhole principle again, $W$ must contain at least
$\lceil \frac{144m(k-1)-143}{12}\rceil$ $=$ $12m(k-1)-11$ $S$-points such that no $I$-point came between them (in the order of $P$).
Denote the set of these $S$-points by $J$ and the $I$-points that came before them by $I'$.
Cover the points of the plane that are not NE from any $I'$-point with $|I'|+1\le 12$ wedges\footnote{Note that here we could win about a factor of $2$ by summing up for various values of $|I'|$ instead of bounding it with $11$.} (two of them are actually half-planes, see Figure \ref{fig_important}).
Notice that none of the $S$-points coming after $I'$ can be NE from an $I'$-point, in particular all $J$-points belong to at least one of these wedges. Thus at the time just before adding the next $I$-point (when all the $J$-points are already present), one of these wedges contains at least $\lceil \frac{|J|}{12}\rceil\ge m(k-1)$ $J$-points.
Since these points were all colored with a color which is at most $k-2$, this contradicts the correctness of the algorithm for $k-1$ colors that we applied for the set $S$. This finishes the proof.

\section{Remarks and open problems}
It has been conjectured that if $P$ is cover-decomposable ($m_P(2)$ exists), then it is also cover-decomposable into many coverings ($m_P(k)$ exists).
This paper is yet another evidence that this conjecture holds
Another conjecture is that in fact $m(2)=O(m(k))$, the next step is to verify this for octants.
For more related questions, see the recent surveys \cite{P10} and \cite{PPT11}.

The following claim follows from the paper of Cardinal and Korman \cite{korman}:

\begin{claim}\label{4color} Any finite set of points in $\R^3$ can be colored with $4$ colors such that any translate of a given octant with at least $2$ points contains two different colors and $4$ colors are sometimes needed.
\end{claim}

In fact they prove the stronger statement that $4$ colors are enough even if we change the octant in the statement to any cone. The above statement can be proved easily with our approach as well, for completeness we sketch this proof.

\begin{proof}
Project the points again on the $z=0$ plane and consider them as an ordered point set. We define a graph $G$ and its drawing in the plane as follows. We connect two points if there exists a wedge $W$ that contains exactly these two points at some time $t$. If the two points are in SW-NE position then we connect them by a straight segment, if they are in NW-SE position then by a reverse-$L$ shape (a polygonal line with two segments, one going very close to the top side of the rectangle defined by these two points and the other segment going very close to the right side of the rectangle). This can be done such that there are no intersections between adjacent edges and it is easy to check that there will be no intersections between non-adjacent edges. Thus $G$ will be a planar graph and $4$-coloring the vertices of $G$ gives the needed coloring. 
\end{proof}

This claim and our earlier result about the cover-decomposability of octants raise the following problem.

\begin{problem}
What is the smallest $k$ such that any finite set of points in $\R^3$ can be colored with $3$ colors such that any translate of a given octant with at least $k$ points contains two different colors?
\end{problem}

We know that $k$ is at least $2$ trivially and at most $12$ as for $k=12$ already $2$ colors suffice \cite{KP}.
In case of the related problem where instead of translates of a given octant we consider discs in the plane, even the existence of such a $k$ is unknown \cite{wcfuj,wcf}.

\end{document}